\documentclass{article}

\newcommand{\N}{I\!\! N}

\begin{document}

Tsemo Aristide

C.P. 79067, Gatineau Canada

 J8Y 6V2

tsemo58@yahoo.ca

\bigskip
\bigskip

\centerline{\bf Quadratic categories, Koszul resolutions and
operads.}

\bigskip
\bigskip

\centerline{\bf Abstract.}

The category of quadratic algebras has been endowed by Manin with
two tensor products. These products have been generalized to
quadratic operads by Ginzburg and Kapranov, and to $n$-homogeneous
algebras by Berger. The purpose of this paper is to define an
abstract notion of quadratic category such that the categories of
quadratic algebras and quadratic operads are examples of this
notion. We define Koszul complexes in this setting, and
representations of quadratic categories in the category of
quadratic algebras, and Tannaka quadratic category.

\bigskip
\bigskip

\centerline{\bf Introduction.}

\bigskip
\bigskip

 A quadratic algebra is a quotient of a tensor algebra
$T(V)$ by an ideal $C$ generated by a subspace of $V\otimes V$.
These algebras are ${\N}$-graded algebras generated by their
elements of degree $1$ which satisfy quadratic relations.
Quadratic algebras appear in different domains of mathematics, in
topology, the notion of  Steenrood algebras is defined, in
differential geometry the Clifford algebras are one of the main
tools to study $Spin$-geometry, in group theory, symmetric and
exterior algebra are very useful. The notion of scheme which is a
fundamental object in algebraic geometry is a projective spectrum
of a quadratic algebra.

Cohomology theories are defined in the general context of abelian
category with enough injective objects by applying the $Hom$ to
resolutions. To compute cohomology groups, we need to define
complexes which represent these   resolutions, like the Chevalley
complex in group theory, the Koszul complex in Lie theory, the Bar
complex in associative algebra... Even at this stage, these
canonical complexes are hardly tractable. This has motivated
Priddy in [7] to defined a generalized Koszul complex for
quadratic algebras which allows to compute their cohomology when
it is a resolution. The complex defined by Priddy is a
generalization of the classical Koszul complex defined with the
exterior and symmetric algebras. This is useful in practice since
the Koszul complex is simpler than the Bar resolution.

 The automorphism group of quadratic algebras has been used in
 theoretical physics in the inverse scattering problem, and in low
 dimensional topology. It is in this context that
  Manin has endowed in [6] the category of quadratic algebras
with two tensor products $o$ and $\bullet$, he has also defined
the notion of quadratic dual which represents the Yoneda algebra
of quadratic Koszul algebras. He has shown that there exists
internal $\underline{Hom}$ object in the category of quadratic
algebras endowed with the tensor product $\bullet$.

An operad is an object which encodes operations. These objects
have been defined in homotopy theories, and  in [5] to study
algebraic structures. The Manin tensor products have been adapted
to the theory of quadratic operads by Ginsburg and Kapranov [5],
in their paper they have defined the notion of Koszul resolution
of a quadratic operad, which is an application of the general
Koszul duality defined by Beilinson Ginsburg and Schechmann [1].
Recently, Berger [3] has defined the category of $n$-Koszul
algebras. He has endowed this category with the similar tensor
products defined by Manin.

The purpose of this paper is to defined a general notion of
quadratic category endowed with two tensor products and a duality
which satisfy compatibility conditions.  The categories of
quadratic algebras, quadratic operads, and $n$-Koszul algebras are
examples of quadratic categories. In this context, we show the
following result:

\bigskip

{\bf Theorem.}

{\it Let $(C,o,\bullet,!)$ be a quadratic category, then for each
objects $U$ and $V$ of $C$, $V o U^!$ is an internal $Hom$ of the
tensor category $(C,\bullet)$.}

\bigskip

 We also define the notion of $n$-quadratic operad and show
how we can use it to encode coherence relations for
$n$-categories. Tensor categories has been studied in Saavedra,
Deligne to determine properties of cohomology ring of algebraic
varieties. These authors have defined the notion of Tanaka
category which is an Abelian rigid tensor category endowed with an
exact faithful functor to a category of vector spaces, and have
shown that a Tanaka category is equivalent to the category of
representations of an affine group scheme. We adapt the study of
these authors to the quadratic context by defining a quadratic
Tanaka category to be an Abelian category endowed with an exact
faithful functor to the category of quadratic algebras, we show:

\bigskip

{\bf Theorem.}

{\it
  A quadratic Tanaka category is equivalent to the
category of quadratic representations of an affine group scheme.}

\bigskip

A projective scheme $U$,is a projective spectrum of a quadratic
algebra $T_U$, we can associated to $U$ the Tanaka category $T'_U$
generated by $T_U$. The property of the group $H_U$ whose category
of quadratic representations is equivalent to $T'_U$ seems to by
an interesting object to study.  Many authors have try to define
non commutative algebraic geometry. The algebraic geometry of
sensor category is studied by Deligne  perhaps quadratic Tanaka
categories represent the good framework to noncommutative
algebraic geometry, and the category of quadratic algebras the
motivic category in non commutative algebraic geometry.

\bigskip
\bigskip

{\bf I. Quadratic categories.}

\bigskip

The purpose of this part is to present the notion and properties
of quadratic categories.

\bigskip

 {\bf Definition 1.}

A quadratic category  $(C,o,\bullet)$ is a category $C$ endowed
with the tensors products $o$ and $\bullet$, whose neutral
elements are respectively $I_o$ and $I_{\bullet}$, and whose
associative constraints are respectively $c_o$ and $c_{\bullet}$.
We suppose that the following properties are satisfied: For each
object $U$ of $C$, there exists an object $U^!$ of $C$, morphisms

$$
c_U:I_{\bullet}\rightarrow U o U^!
$$

$$
d_U:U^!\bullet U\rightarrow I_0
$$

$$
f_{U_1,U_2,U_3}: (U_1 o U_2)\bullet U_3\longrightarrow U_1 o (U_2
\bullet U_3)
$$

$$
h_{U_1,U_2,U_3}:U_1\bullet (U_2 o U_3)\longrightarrow (U_1 \bullet
U_2) o U_3
$$

such that the following diagrams are commutative:

$$
\matrix{(U_1\bullet (U_2 o U_3))\bullet U_4 &{\buildrel
{c_{\bullet}(U_1,U_2 o U_3, U_4)}\over{\longrightarrow}}
U_1\bullet ((U_2 o U_3)\bullet U_4){\buildrel {Id_{U_1}\bullet
f_{U_2,U_3,U_4}}\over {\longrightarrow}} & U_1\bullet (U_2 o
(U_3\bullet U_4))\cr \downarrow h_{U_1, U_2, U_3}\bullet
Id_{U_4}&& \downarrow h_{U_1,U_2, U_3\bullet U_4} \cr ((U_1\bullet
U_2) o U_3)\bullet U_4  &{\buildrel{f_{U_1\bullet
U_2,U_3,U_4}}\over{\longrightarrow}}& (U_1\bullet U_2) o
(U_3\bullet U_4)}\leqno (1)
$$

$$
\matrix{(U_1 o U_2)\bullet (U_3 o U_4)&{\buildrel{f_{U_1,U_2,U_3 o
U_4}}\over{\longrightarrow}} &U_1 o (U_2\bullet(U_3 o U_4))\cr
\downarrow h_{U_1 o U_2,U_3,U_4}&&\downarrow Id_{U_1} o
h_{U_2,U_3,U_4}\cr ((U_1 o U_2)\bullet U_3) o U_4
&{\buildrel{f_{U_1,U_2,U_3}o Id_{U_4}}\over{\longrightarrow}} (U_1
o(U_2\bullet U_3)) o U_4 {\buildrel{c_0(U_1,U_2\bullet
U_3,U_4)}\over{\longrightarrow}} & U_1 o ((U_2\bullet U_3) o
U_4)}\leqno (2)
$$

$$
\matrix{I_{\bullet}\bullet U&{\buildrel{c_U\bullet
Id_U}\over{\longrightarrow}} (U o U^!)\bullet
U{\buildrel{f_{U,U^!,U}}\over{\longrightarrow}} U o (U^!\bullet
U){\buildrel{Id_U o d_U}\over{\longrightarrow}}& U o I_o\cr
\downarrow &&\downarrow \cr
U&{\buildrel{Id_U}\over{\longrightarrow}} &U} \leqno (3)
$$

$$
\matrix{U^!\bullet I_{\bullet}&{\buildrel{Id_{U^!}\bullet
c_U}\over{\longrightarrow}} U^! \bullet (U o
U^!){\buildrel{h_{U^!,U,U^!}}\over{\longrightarrow}} (U^!\bullet
U) o U^!{\buildrel{{d_U} o Id_{U^!}}\over{\longrightarrow}}& I_o o
U^!\cr \downarrow &&\downarrow \cr
U^!&{\buildrel{Id_U}\over{\longrightarrow}} &U^!}\leqno (4)
$$

Let $u_1:U_1\rightarrow U'_1$, $u_2:U_2\rightarrow U'_2$ and
$u_3:U_3\rightarrow U'_3$ be three morphisms of $C$. The following
diagrams are supposed to be commutative:

$$
\matrix{U_1\bullet (U_2 o U_3)&{\buildrel{u_1\bullet (u_2 o
u_3)}\over{\longrightarrow}}& U'_1\bullet (U'_2 o U'_3)\cr
\downarrow h_{U_1,U_2,U_3}&&\downarrow h_{U'_1,U'_2,U'_3}\cr
(U_1\bullet U_2) o U_3 &{\buildrel{(u_1\bullet u_2)o
u_3}\over{\longrightarrow}}&(U'_1\bullet U'_2) o U'_3}\leqno (5)
$$

$$
\matrix{(U_1 o U_2)\bullet U_3 &{\buildrel{(u_1 o u_2)\bullet
u_3}\over{\longrightarrow}} &(U'_1 o U'_2)\bullet
U'_3\cr\downarrow f_{U_1,U_2,U_3}&&\downarrow
f_{U'_1,U'_2,U'_3}\cr U_1 o (U_2\bullet U_3)&{\buildrel{u_1 o
(u_2\bullet u_3)}\over{\longrightarrow}}& (U'_1 o U'_2)\bullet
U'_3}\leqno (6)
$$

The example who has motivated the construction of quadratic
categories is the theory of quadratic algebras, recall the
constructions of the tensor products defined by Manin.

Let $U=T(U_1)/(C_1)$ and $V=T(V_1)/(C_2)$ be two quadratic
algebras respectively isomorphic to the quotient of the tensor
algebra of the finite dimensional $L$-vector spaces $U_1$ and
$V_1$ by the ideals $C_1\subset U_1\otimes U_1$ and $C_2\subset
V_1\otimes V_1$. We endow the class of quadratic algebra with the
structure of a category such that $Hom(U,V)$ is the set of
morphisms of algebras $h:U\rightarrow V$ defined by a linear
application $h_1:U_1\rightarrow V_1$ such that $(h_1\otimes
h_1)(C_1)\subset C_2$.

We define $U\bullet V$ to be the quadratic algebra isomorphic to
the quotient of the tensor algebra of $T(U_1\otimes V_1)$ by the
ideal generated by $t_{23}(C_1\otimes C_2)$. The isomorphism
$t_{23}$ of $U_1^{\otimes^2}\otimes V_1^{\otimes^2}$ is defined by
$t_{23}(u_1\otimes u'_1\otimes v_1\otimes v'_1)=u_1\otimes
v_1\otimes u'_1\otimes v'_1$.

The tensor product $U o V$ is defined to be the quotient of the
tensor algebra of $T(U_1\otimes V_1)$ by the ideal generated by
$t_{23}(U_1^{\otimes^2}\otimes C_2+  C_1\otimes V_1^{\otimes^2})$.

\medskip

The quadratic dual $U^!$ of $U$ is the quadratic algebra
$T(U_1^*)/C_1^*$ where $U_1^*$ is the dual vector space of $U_1$,
and $C_1^*$ the annihilator of $C_1$ in $U^*_1\otimes U^*_1$.

\medskip

The neutral element $I_o$ is $T(L)$, and $I_{\bullet}$ is $L$. The
map $c_U:I_{\bullet}\rightarrow U o U^!$ is defined $L\rightarrow
U_1\otimes U_1^*$ $1\rightarrow \sum_{i=1}^{i=n}u_i\otimes u^i$,
where $(u_1,..,u_n)$ is a basis of $U_1$, $U_1^*$ the dual space
of $U_1$ and $(u^1,..,u^n)$ its dual basis. The map
$d_U:U^!\bullet U\rightarrow I_o$ is defined by the duality
$U_1\otimes U_1^*\rightarrow L$.

\medskip

Let $U^i=T(U^i)/(C_i), i=1,2,3$. The ideal $C$ which defines $(U^1
o U^2)\bullet U^3$ is $t_{23}(t_{23}(U_1^{\otimes^2}\otimes C_2 +
C_1\otimes U_2^{\otimes^2})\otimes C_3)$, and the ideal $C'$ which
defines $U^1 o (U^2\bullet U^3)$ is
$t_{23}(U_1^{\otimes^2}\otimes(t_{23}(C_2\otimes C_3))+C_1\otimes
(U_2\otimes U_3)^{\otimes^2}$). We remark that $C$ is contained in
$C'$, the associative constraint $(U_1\otimes U_2)\otimes
U_3\longrightarrow U_1\otimes (U_2\otimes U_3)$ of vector spaces
defines the quadratic constraint $f_{U_1,U_2,U_3}$.

The ideal $D$ which defines $U^1\bullet (U^2 o U^3)$ is
$t_{23}(C_1\otimes t_{23}(U_2^{\otimes^2}\otimes C_3 + C_2\otimes
U_3^{\otimes^2}))$. The ideal $D'$ which defines the algebra
$(U^1\bullet U^2) o U^3$ is $t_{23}(({U_1\otimes
U_2}^{\otimes^2}\otimes C_3+t_{23}(C_1\otimes C_2)\otimes
U_3^{\otimes^2})$. We remark that $D'$ contains $D$ this implies
that the associative constraint for vector spaces $U_1\otimes
(U_2\otimes U_3)\longrightarrow (U_1\otimes U_2)\otimes U_3$
defined the quadratic associative constraint $h_{U_1,U_2,U_3}$.

\medskip

\bigskip

{\bf Theorem 2.}

{\it Let $C$ be a quadratic category, and $L$ and $N$ two objects
of $C$. The opposite functor:

$$
C\longrightarrow C
$$

$$
U\longrightarrow Hom(U\bullet L,N)
$$

is representable by $N o {L}^!$.}

\bigskip

{\bf Proof.}

We define a morphism between the functors $U\rightarrow
Hom(U\bullet L,N)$ and $U\rightarrow Hom(U,N o L^!)$ by defining
for each element $u$ of $Hom(U\bullet L,N)$, the element $u'$ in
$Hom(U,N o L^!)$ by:

$$
\matrix{U\rightarrow U\bullet I_{\bullet}{\buildrel{Id_U\bullet
c_L}\over{\rightarrow}} U\bullet(L o
L^!){\buildrel{h_{U,L,L^!}}\over {\longrightarrow}} (U\bullet L) o
L^!{\buildrel{u o Id_{L^!}}\over{\longrightarrow}} N o L^!}
$$

We define a morphism between the functors $U\rightarrow Hom(U,N o
L^!)$ and $U\rightarrow Hom(U\bullet N,L)$ by assigning to each
element $v$ in $Hom(U,N o L^!)$ the element $v"$ in $Hom(U\bullet
L,N)$ defined by:

$$
\matrix{U\bullet L {\buildrel{v\bullet
Id_L}\over{\longrightarrow}}(N o L^!)\bullet
L{\buildrel{f_{N,L^!,L}}\over{\longrightarrow}} N o (L^!\bullet
L){\buildrel{Id_N o d_L}\over{\longrightarrow}}}N
$$

We have to show that the correspondence defined on $Hom(U\bullet
L,N)$ by $u\rightarrow (u')"$ is  the identity of the functor
$U\rightarrow Hom(U\bullet L,N)$. We have:

$$
(u')"=\matrix{(U\rightarrow U\bullet
I_{\bullet}{\buildrel{Id_U\bullet c_L}\over{\rightarrow}}
U\bullet(L o L^!){\buildrel{h_{U,L,L^!}}\over {\longrightarrow}}
(U\bullet L) o L^!{\buildrel{u o Id_{L^!}}\over{\longrightarrow}}
N o L^!)\bullet L}
$$

$$
\matrix{{\buildrel{f_{N,L^!,L}}\over{\longrightarrow}} N o
(L^!\bullet L){\buildrel{Id_N o d_L}\over{\longrightarrow}}}N
$$

Let $(C,\otimes)$ be a tensor category, and $u:U\rightarrow U'$,
$u':U'\rightarrow U"$, $v:V\rightarrow V'$, $v":V'\rightarrow V"$
arrows of $C$, $(u'\otimes v')\circ (u\otimes v)=(u'\circ
u)\otimes (v'\circ v)$. Applying this fact to $\bullet$ at the
first line of the previous equality, we obtain:

$$
(u')"=\matrix{U\bullet L\rightarrow (U\bullet I_{\bullet})\bullet
L{\buildrel{Id_U\bullet c_L\bullet Id_L}\over{\rightarrow}}
(U\bullet(L o L^!))\bullet L{\buildrel{h_{U,L,L^!}\bullet
Id_L}\over {\longrightarrow}} ((U\bullet L) o L^!)\bullet
L{\buildrel{u o Id_{L^!}\bullet Id_L}\over{\longrightarrow}} (N o
L^!)\bullet L}
$$

$$
\matrix{{\buildrel{f_{N,L^!,L}}\over{\longrightarrow}} N o
(L^!\bullet L){\buildrel{Id_N o d_L}\over{\longrightarrow}}N}
$$

 Applying property $(5)$ we obtain
 $$
 \matrix{((U\bullet L) o L^!)\bullet
L{\buildrel{u o Id_{L^!}\bullet Id_L}\over{\longrightarrow}} (N o
L^!)\bullet L {\buildrel{f_{N,L^!,L}}\over{\longrightarrow}} N o
(L^!\bullet L){\buildrel{Id_N o d_L}\over{\longrightarrow}}N}
$$

 $$
 =\matrix{((U\bullet L) o L^!)\bullet
L \matrix{{\buildrel{f_{U\bullet L,  L^!,L}}\over
{\longrightarrow}} (U\bullet L) o (L^!\bullet L){\buildrel {u o
Id_{L^!\bullet L}}\over {\longrightarrow}}
 N o (L^!\bullet L){\buildrel{Id_N o d_L }\over{\longrightarrow}}N}}
$$

We deduce that
$$
(u')"=\matrix{U\bullet L\rightarrow (U\bullet I_{\bullet})\bullet
L{\buildrel{Id_U\bullet c_L\bullet Id_L}\over{\rightarrow}}
(U\bullet(L o L^!))\bullet L{\buildrel{h_{U,L,L^!}\bullet
Id_L}\over {\longrightarrow}} ((U\bullet L) o L^!)\bullet L}
$$

$$
\matrix{{\buildrel{f_{U\bullet L,  L^!,L}}\over {\longrightarrow}}
(U\bullet L) o (L^!\bullet L){\buildrel {u o Id_{L^!\bullet
L}}\over {\longrightarrow}}
 N o (L^!\bullet L){\buildrel{Id_N o d_L }\over{\longrightarrow}}N}
$$
 Applying $(1)$ we obtain that

$$
\matrix{(U\bullet(L o L^!))\bullet L{\buildrel{h_{U,L,L^!}\bullet
Id_L}\over {\longrightarrow}} ((U\bullet L) o L^!)\bullet L
{\buildrel{f_{U\bullet L,  L^!,L}}\over {\longrightarrow}}
(U\bullet L) o (L^!\bullet L)}
$$

$$
=\matrix{(U\bullet(L o L^!))\bullet L{\buildrel{{c_{\bullet}}{U,L
o L^!,L}}\over{\longrightarrow}} U\bullet ((L o L^!)\bullet
L){\buildrel{Id_U o f_{L,L^!,L}}\over {\longrightarrow}}U\bullet(L
o (L^!\bullet L)){\buildrel{h_{U,L,L^!\bullet L}}\over
{\longrightarrow}} (U\bullet L) o (L^!\bullet L)}
$$

We deduce that

$$
(u')"=\matrix{U\bullet L\rightarrow (U\bullet I_{\bullet})\bullet
L{\buildrel{Id_U\bullet c_L\bullet Id_L}\over{\rightarrow}}
(U\bullet(L o L^!))\bullet L{\buildrel{{c_{\bullet}}_{U,L o
L^!,L}}\over{\longrightarrow}} U\bullet ((L o L^!)\bullet L)}
$$

$$
\matrix{{\buildrel{Id_U o f_{L,L^!,L}}\over
{\longrightarrow}}U\bullet(L o (L^!\bullet
L)){\buildrel{h_{U,L,L^!\bullet L}}\over {\longrightarrow}}
(U\bullet L) o (L^!\bullet L){\buildrel{u o Id_{L^!\bullet
L}}\over {\longrightarrow}} N o (L^!\bullet L){\buildrel {Id_N o
d_L}\over{\longrightarrow}}N }
$$

Using the fact that $o$ is a tensor product we deduce that:

$$
(u')"=\matrix{U\bullet L\rightarrow (U\bullet I_{\bullet})\bullet
L{\buildrel{Id_U\bullet c_L\bullet Id_L}\over{\rightarrow}}
(U\bullet(L o L^!))\bullet L{\buildrel{{c_{\bullet}}_{U,L o
L^!,L}}\over{\longrightarrow}} U\bullet ((L o L^!)\bullet L)}
$$

$$
\matrix{{\buildrel{Id_U o f_{L,L^!,L}}\over
{\longrightarrow}}U\bullet(L o (L^!\bullet
L)){\buildrel{h_{U,L,L^!\bullet L}}\over {\longrightarrow}}
(U\bullet L) o (L^!\bullet L){\buildrel{Id_L o d_L}\over
{\longrightarrow}} (U\bullet L) o I_o {\buildrel
u\over{\longrightarrow}}N }
$$
Using property $(5)$, we obtain that:

$$
\matrix{U\bullet(L o (L^!\bullet L)){\buildrel{h_{U,L,L^!\bullet
L}}\over {\longrightarrow}} (U\bullet L) o (L^!\bullet
L){\buildrel{Id_L o d_L}\over {\longrightarrow}} (U\bullet L) o
I_o}
$$

$$
=\matrix{U\bullet (L o (L^!\bullet L)){\buildrel{Id_U\bullet(Id_L
o d_L)}\over{\longrightarrow}} U\bullet (L o I_o){\buildrel
{h_{U,L,I_o}}\over{\longrightarrow}}(U\bullet L) o I_o}
$$

We deduce that

$$
(u')"=\matrix{U\bullet L\rightarrow (U\bullet I_{\bullet})\bullet
L{\buildrel{Id_U\bullet c_L\bullet Id_L}\over{\rightarrow}}
(U\bullet(L o L^!))\bullet L{\buildrel{{c_{\bullet}}_{U,L o
L^!,L}}\over{\longrightarrow}} U\bullet ((L o L^!)\bullet L)}
$$

$$
\matrix{{\buildrel{Id_U o f_{L,L^!,L}}\over
{\longrightarrow}}U\bullet (L o (L^!\bullet
L)){\buildrel{Id_U\bullet(Id_L o d_L)}\over{\longrightarrow}}
U\bullet (L o I_o){\buildrel
{h_{U,L,I_o}}\over{\longrightarrow}}(U\bullet L) o I_o{\buildrel
u\over{\longrightarrow}}N }
$$

Applying $(5)$, we obtain that

$$
\matrix{ ((U\bullet I_{\bullet})\bullet L{\buildrel{Id_U\bullet
c_L\bullet Id_L}\over{\rightarrow}} (U\bullet(L o L^!))\bullet
L{\buildrel{{c_{\bullet}}_{U,L o L^!,L}}\over{\longrightarrow}}
U\bullet ((L o L^!)\bullet L)}
$$

$$
=\matrix{(U\bullet I_{\bullet})\bullet
L{\buildrel{{c_{\bullet}}{U,I_{\bullet},L}}\over{\longrightarrow}}U\bullet
(I_{\bullet}\bullet L){\buildrel{Id_U\bullet (c_L\bullet
Id_L)}\over{\longrightarrow}} U\bullet((L o L^!)\bullet L)}
$$
We deduce that

$$
(u')"=\matrix{U\bullet L\rightarrow (U\bullet I_{\bullet})\bullet
L {\buildrel{h_{U,I_{\bullet},L}}\over{\longrightarrow}}U\bullet
(I_{\bullet}\bullet L){\buildrel{Id_U\bullet (c_L\bullet
Id_L)}\over{\longrightarrow}} U\bullet((L o L^!)\bullet L)}
$$

$$
\matrix{{\buildrel{Id_U\bullet{h_{L,L^!,L}}}\over{\longrightarrow}}U\bullet
(L o (L^!\bullet L)){\buildrel{Id_U \bullet (Id_L o
d_L)}\over{\longrightarrow}} U\bullet (L o I_o){\buildrel
{h_{U,L,I_o}}\over{\longrightarrow}}(U\bullet L) o
I_o{\buildrel{u}\over{\longrightarrow}} N}
$$
 Applying $(3)$ we deduce
 that

 $$
 \matrix{U \bullet L\rightarrow U\bullet
(I_{\bullet}\bullet L){\buildrel{Id_U\bullet (c_L\bullet
Id_L)}\over{\longrightarrow}} U\bullet((L o L^!)\bullet L)}
$$

$$
\matrix{
{\buildrel{Id_U\bullet{h_{L,L^!,L}}}\over{\longrightarrow}}U\bullet
(L o (L^!\bullet L)){\buildrel{Id_U \bullet (Id_L o
d_L)}\over{\longrightarrow}} U\bullet (L o I_o)\longrightarrow
U\bullet L}
$$

$$
= Id_{U\bullet L}: U\bullet  L\longrightarrow U\bullet L
$$

This implies that:

$$
(u')"=\matrix{U\bullet L{\buildrel {u}\over{\longrightarrow}} N}
$$

We have to show now that $(v")'=v$. We have:

$$
(v")'=\matrix{U\longrightarrow U\bullet
I_{\bullet}{\buildrel{Id_U\bullet c_L}\over{\longrightarrow}}
U\bullet (L o L^!){\buildrel
{h_{U,L,L^!}}\over{\longrightarrow}}}(U\bullet L) o L^!
$$

$$
((U\bullet L){\buildrel {v\bullet Id_L}\over{\longrightarrow}} (N
o L^!)\bullet L {\buildrel{f_{N,L^!,L}}\over{\longrightarrow}} N o
(L^!\bullet L){\buildrel{Id_N\bullet d_L}\over{\longrightarrow
N}})o L^!\longrightarrow N o L^!
$$

Using the fact that $o$ is a tensor product, we obtain that:

$$
\matrix{((U\bullet L){\buildrel {u\bullet
Id_L}\over{\longrightarrow}} (N o L^!)\bullet L
{\buildrel{f_{N,L^!,L}}\over{\longrightarrow}} N o (L^!\bullet
L){\buildrel{Id_N\bullet d_L}\over{\longrightarrow N}})o L^!}
$$

$$
=\matrix{((U\bullet L) o L^!{\buildrel {(v\bullet Id_L) o
Id_{L^!}}\over{\longrightarrow}} ((N o L^!)\bullet L) o L^!
{\buildrel{f_{N,L^!,L} o Id_{L^!}}\over{\longrightarrow}} (N o
(L^!\bullet L)) o L^!{\buildrel{(Id_N o d_L) o
Id_{L^!}}\over{\longrightarrow}} No L^!}
$$

This implies that

$$
(v")'=\matrix{U\longrightarrow U\bullet
I_{\bullet}{\buildrel{Id_U\bullet c_L}\over{\longrightarrow}}
U\bullet (L o L^!){\buildrel {h_{U,L,L^!}}\over{\longrightarrow}}}
$$

$$
\matrix{(U\bullet L) o L^!{\buildrel {v o
Id_{L^!}}\over{\longrightarrow}} ((N o L^!)\bullet L) o L^!
{\buildrel{f_{N,L^!,L} o Id_{L^!}}\over{\longrightarrow}} (N o
(L^!\bullet L)) o L^!{\buildrel{(Id_N o d_L) o Id_{L^!}}\over
{\longrightarrow}} No L^!}
$$

Using property $(5)$ and the fact that $\bullet$, $o$ are tensor
products, we obtain that:

$$
\matrix{U\longrightarrow U\bullet
I_{\bullet}{\buildrel{Id_U\bullet c_L}\over{\longrightarrow}}
U\bullet (L o L^!){\buildrel
{h_{U,L,L^!}}\over{\longrightarrow}}(U\bullet L) o L^!{\buildrel
{(v\bullet Id_L) o Id_{L^!}}\over{\longrightarrow}} ((N o
L^!)\bullet L) o L^!}
$$

$$
=\matrix{U{\buildrel v\over {\longrightarrow}} (N o L^!) \bullet
I_{\bullet}{\buildrel{Id_{N o L^!}\bullet
c_L}\over{\longrightarrow}} (N o L^!)\bullet (L o L^!){\buildrel
{h_{N o L^!,L,L^!}}\over{\longrightarrow}}((N o L^!)\bullet L) o
L^!}
$$

We deduce that

$$
(v")'=\matrix{U{\buildrel v\over {\longrightarrow}} (N o L^!)
\bullet I_{\bullet}{\buildrel{Id_{N o L^!}\bullet
c_L}\over{\longrightarrow}} (N o L^!)\bullet (L o L^!){\buildrel
{h_{N o L^!,L,L^!}}\over{\longrightarrow}}((N o L^!)\bullet L) o
L^!}
$$

$$
\matrix{{\buildrel{f_{N,L^!,L} o Id_{L^!}}\over{\longrightarrow}}
(N o (L^!\bullet L)) o L^!{\buildrel{{c_o}(N,L^!\bullet
L,L^!)}\over{\longrightarrow}}N o ((L^!\bullet L)o
L^!)\longrightarrow No L^!}
$$

The property $(2)$ implies that:

$$
\matrix{(N o L^!)\bullet (L o L^!){\buildrel{h_{N o
L^!,L,L^!}}\over{\longrightarrow}}((N o L^!)\bullet L) o L^!
{\buildrel{f_{N,L^!,L} o Id_{L^!}}\over{\longrightarrow}} (N o
(L^!\bullet L)) o L^! {\buildrel{{c_o}_{N,L^!\bullet
L,L^!}}\over{\longrightarrow}}N o ((L\bullet L^!) o L)}
$$

$$
=\matrix{ (N o L^!)\bullet (L o L^!) {\buildrel{f_{N,L^!,L o
L^!}}\over{\longrightarrow}} N o (L^!\bullet (L o
L^!)){\buildrel{Id_N \bullet
h_{L^!,L,L^!}}\over{\longrightarrow}}N o ((L^!\bullet L) o L)}
$$

This implies that:

$$
(v")'=\matrix{U{\buildrel v\over {\longrightarrow}} (N o L^!)
\bullet I_{\bullet}{\buildrel{Id_{N o L^!}\bullet
c_L}\over{\longrightarrow}} (N o L^!)\bullet (L o L^!)}
$$

$$
\matrix{{\buildrel{f_{N,L^!,L o L^!}}\over{\longrightarrow}} N o
(L^!\bullet (L o L^!)){\buildrel{Id_N \bullet
h_{L^!,L,L^!}}\over{\longrightarrow}}N o ((L^!\bullet L) o
L)\rightarrow N oL^!}
$$

The property $(5)$implies that:

$$
\matrix{(N o L^!) \bullet I_{\bullet}{\buildrel{Id_{N o
L^!}\bullet c_L}\over{\longrightarrow}} (N o L^!)\bullet (L o L^!)
{\buildrel{f_{N,L^!,L o L^!}}\over{\longrightarrow}} N o
(L^!\bullet (L o L^!))}
$$

$$
=\matrix{(N o L^!) \bullet
I_{\bullet}{\buildrel{f_{N,L,I_{\bullet}}}\over{\longrightarrow}}
N  o (L^!\bullet I_{\bullet}){\buildrel{Id_N o(Id_{L^!}\bullet
c_L)}\over{\longrightarrow}} N o (L ^!\bullet(L o L^!))}
$$

This implies that:

$$
(v")'=\matrix{U{\buildrel v\over {\longrightarrow}}(N o L^!)
\bullet
I_{\bullet}{\buildrel{f_{N,L,I_{\bullet}}}\over{\longrightarrow}}
N  o (L^!\bullet I_{\bullet})}
$$

$$
\matrix{{\buildrel{Id_N o(Id_L\bullet c_L)}\over{\longrightarrow}}
N o (L^! \bullet (L o L^!)){\buildrel{Id_N o
h_{L^!,L,L^!}}\over{\longrightarrow}} N o ((L^!\bullet L) o
L^!)\rightarrow N o L^!}
$$

The property $(4)$ implies that:

$$
\matrix{L^!\longrightarrow L^! \bullet
I_{\bullet}{\buildrel{Id_{L^!}\bullet
c_L}\over{\longrightarrow}}L^!\bullet (L o
L^!){\buildrel{h_{L,L^!,L}}\over{\longrightarrow}}(L^!\bullet L) o
L^!\longrightarrow L^!}
$$

is the identity of $L^!$. This implies that $(v")'=v$ $\bullet$

\bigskip

The previous theorem implies that the tensor category
$(C,\bullet)$ is endowed with an internal $\underline{Hom}(U,V)= V
o L^!$. The general properties of internal $Hom$ implies the
existence of an isomorphism

$$
d_{U_1,U_2}:\underline{Hom}(U_1,U_2)\bullet U_1\rightarrow U_2
$$

and a map

$$
l_{U_1,U_2,U_3}:\underline{Hom}(U_2,U_3)\bullet\underline{Hom}(U_1,U_2)
\longrightarrow \underline{Hom}(U_1,U_3)
$$

The map $d_{U_1,U_2}$ is the composition:

$$
\matrix{(U_2 o U_1^!)\bullet
U_1{\buildrel{f_{U_2,U_1,U_1^!}}\over{\longrightarrow}} U_2 o
(U_1^!\bullet U_1){\buildrel{Id_{U_2} o
d_{U_1}}\over{\longrightarrow}} U_2 o I_o\longrightarrow U_2}
$$

The map $l_{U_1,U_2,U_3}$ is the composition:

$$
\matrix{(U_3 o U_2^!)\bullet (U_2 o
U_1^!){\buildrel{f_{U_3,U_2^!,U_2 o U_1^!}}\over{\longrightarrow}}
U_3 o (U_2^!\bullet (U_2 o U_1^!)){\buildrel{Id_{U_3} o
h_{U_2^!,U_2,U_1}}\over{\longrightarrow}} U_3 o ((U_2^!\bullet
U_2) o U_1^!){\buildrel{Id_{U_3} o (d_{U_2} o
Id_{U_1^!})}\over{\longrightarrow}} U_3 o U_1^!}.
$$

\bigskip

{\bf Proposition 3.}

{\it Let $(C,\bullet, o)$ be a quadratic category. There exists a
canonical isomorphism $(U\bullet V)^!\longrightarrow V^! o U^!$,
and the correspondence $U\rightarrow U^!$ is an opposite functor.}

\bigskip

{\bf Proof.}

The general properties of tensor categories imply the existence of
an isomorphism:
$$
\underline{Hom}(U_1\bullet
U_2,U_3)\longrightarrow\underline{Hom}(U_1,\underline{Hom}(U_2,U_3))
$$

we obtain a  isomorphism:

$$
U_3 o (U_1\bullet U_2)^!\longrightarrow (U_3 o U_2^!) o U_1^!
$$

Suppose that $U_3=I_o$  we deduce the existence of an isomorphism:

$$
(U_1\bullet U_2)^!\longrightarrow U_2^! o U_1^!
$$

Let $h:U\rightarrow U'$ be a map. We define the morphism
$h^!:{U'}^!\rightarrow U^!$ as follows:

$$
\matrix{{U'}^!\rightarrow {U'}^!\bullet
I_{\bullet}{\buildrel{Id_{{U'}^!}\bullet
c_U}\over{\longrightarrow}}{U'}^!\bullet (U o U^!){\buildrel
{h_{{U'}^!,U,U^!}}\over{\longrightarrow}}({U'}^!\bullet U) o
U^!{\buildrel{(Id_{{U'}^!}\bullet h)o
Id_{U^!}}\over{\longrightarrow}} ({U'}^!\bullet U') o
U^!{\buildrel{d_{U'} o I_{U^!}}\over{\longrightarrow}} U^!}
$$
\medskip

{\bf Proposition 4.}

{\it Let $(C,o,\bullet)$ be a quadratic category, for every object
$U$ of $C$, $U^!$ is unique up to an isomorphism.}

\medskip

{\bf proof.}

The object $U^!$ represents the opposite functor $V\rightarrow
Hom_C(V\bullet U,I_o)$. This implies that this object is unique up
to an isomorphism.

\bigskip

We denote by $\underline{hom}(U,V)$ the object
$\underline{Hom}(U,V)$. We apply the previous proposition to the
map $l_{U_1,U_2,U_3}$. We obtain a map:

$$
l^!_{U_1,U_2,U_3}:\underline{hom}(U_1,U_3)\longrightarrow
\underline{hom}(U_1,U_2) o \underline{hom}(U_2,U_3)
$$

The morphism $l_{U,U,U}$ endows $\underline{Hom}(U,U)$ with a
product
$l_{U,U,U}=l_U:\underline{Hom}(U,U)\bullet\underline{Hom}(U,U)\rightarrow
\underline{Hom}(U,U)$, and the morphism $l^!_{U,U,U}={l^!}_U:
\underline{hom}(U,U)\longrightarrow \underline{hom}(U,U)
o\underline{hom}(U,U)$ endows $\underline{hom}(U,U)$ with a
coproduct. There exists a canonical product on
$\underline{hom}(U,U)$ defined by

$$
\matrix{(U^!\bullet U) o (U^!\bullet U){\buildrel{Id_{U^!\bullet
U} o d_U}\over{\longrightarrow}} U^!\bullet U}
$$

 We can deduce the fact that the correspondence $U\rightarrow U^!$ is
 functorial from general properties of tensor categories since the
 tensor product $\bullet$ has an internal $Hom$ $\bullet$

\bigskip

{\bf Proposition -Definition.}

{\it The object $(\underline{Hom}(U,U), l_U)$ is an algebra, and
the object $(\underline{hom}(U,U), {l^!}_U)$, is a cogebra. This
means that the following diagrams commute:

$$
\matrix{\underline{Hom}(U,U)\bullet
(\underline{Hom}(U,U)\bullet\underline{Hom}(U,U))
{\buildrel{Id_{\underline{Hom}(U,U)}\bullet l_U}\over
{\longrightarrow}}
\underline{Hom}(U,U)\bullet\underline{Hom}(U,U)\cr\cr \downarrow
c_{\bullet}(\underline{Hom}(U,U),\underline{Hom}(U,U),\underline{Hom}(U,U))\downarrow
l_U\cr\cr (\underline{Hom}(U,U)\bullet\underline{Hom}(U,U))\bullet
\underline{Hom}(U,U){\buildrel{l_U\bullet
Id_{\underline{Hom}(U,U)}}\over
{\longrightarrow}}\underline{Hom}(U,U)\bullet\underline{Hom}(U,U)
{\buildrel{l_U}\over{\longrightarrow}} \underline{Hom}(U,U)}
$$
\bigskip

$$
\matrix{\underline{Hom}(U,U)\bullet
I_{\bullet}{\buildrel{Id_{\underline{Hom}(U,U)}\bullet
c_{\underline{Hom}(U,U)}}\over{\longrightarrow}}&\underline{Hom}(U,U)\bullet
\underline{Hom}(U,U)\cr\cr \downarrow &\downarrow l_U\cr\cr
\underline{Hom}(U,U){\buildrel{Id_{\underline{Hom}(U,U)}}\over{\longrightarrow}}
&\underline{Hom}(U,U)}
$$

These two diagrams define the algebra structure of
$\underline{Hom}(U,U)$. The next define the coalgebra structure of
$\underline{hom}(U,U)$.

$$
\matrix{\underline{hom}(U,U){\buildrel{l^!_U}\over{\longrightarrow}}
&\underline{hom}(U,U) o \underline{hom}(U,U)\cr\cr \downarrow l_U
&\downarrow l^!_U o Id_{\underline{hom}(U,U)}\cr\cr
&(\underline{hom}(U,U)\underline{hom}(U,U)) o \underline{hom}(U,U)
\cr &\downarrow c_o\cr\cr\underline{hom}(U,U) o
\underline{hom}(U,U)&{\buildrel{Id_{\underline{hom}(U,U)} o
l^!_U}\over{\longrightarrow}}
\underline{hom}(U,U)(\underline{hom}(U,U) o \underline{hom}(U,U))}
$$

\bigskip

$$
\matrix{\underline{hom}(U,U){\buildrel{l^!_U}\over{\longrightarrow}}
&\underline{hom}(U,U) o \underline{hom}(U,U)\cr\cr \downarrow
&\downarrow d_{\underline{hom}(U,U)} o
d_{\underline{hom}(U,U)}\cr\cr \underline{hom}(U,U)\longrightarrow
&\underline{hom}(U,U)}
$$}

\bigskip

{\bf Proof}

We have only to show that the two first diagrams are commutative
since the two last are their dual. The fact that first two
diagrams commute follows from general properties of tensor
categories $\bullet$

\bigskip

{\bf Proposition 6.}

{\it Let $(C,\bullet,o,!)$ be a quadratic category. The quadratic
dual $I_{\bullet}^!$ is isomorphic to $I_o$, and $I_o^!$ is
isomorphic to $I_{\bullet}$.}

\medskip

{\bf Proof.}

The theorem 2 implies that the functor defined on $C$ by
$U\rightarrow Hom_C(U,I_o)=Hom_C(U\bullet I_{\bullet},I_o)$ is
representable by $I_o\circ I_{\bullet}^!$ and $I_o$. Since
$I_o\circ I_{\bullet}^!$ is isomorphic to $I_{\bullet}^!$, we
deduce that $I_{\bullet}^!$ is isomorphic to $I_o$.

\bigskip

{\bf Definition 7.}

Let $(C, o,\bullet,!)$ be a quadratic category, the commutative
constraints $c'_{\bullet}$ of $(C,\bullet)$ and $c'_o$ of $(C,o)$
is a quadratic braiding if for every objects $U_1$, $U_2$ and
$U_3$ of $C$, the following diagram is commutative:

$$
\matrix{U_1\bullet(U_2 o U_3)&{\buildrel{c'_{\bullet}(U_1,U_2 o
U_3)}\over{\longrightarrow}}(U_2 o U_3)\bullet
U_1{\buildrel{c'_o(U_2,U_3)\bullet
Id_{U_1}}\over{\longrightarrow}}&(U_3 o U_2)\bullet U_1\cr
\downarrow h_{U_1,U_2,U_3}&&\downarrow f_{U_3,U_2,U_1}\cr
(U_1\bullet U_2) o U_3 &{\buildrel{c'_o(U_1\bullet
U_2,U_3)}\over{\longrightarrow}}U_3 o (U_1\bullet
U_2){\buildrel{Id_{U_3} o
c'_{\bullet}(U_1,U_2)}\over{\longrightarrow}}& U_3 o(U_2\bullet
U_1)}
$$

\bigskip

{\bf Proposition 8.}

{\it Let $(C,o,\bullet,!)$ be a quadratic category endowed with a
quadratic braiding, and $U$ an object of $C$, then ${U^!}^!$ is
isomorphic to $U$. In particular ${I^!_{\bullet}}^!\simeq
{I^!_o}\simeq I_{\bullet}$.}

\medskip

{\bf Proof.}

Let $U$ be an object of $C$, consider the morphisms

$$
\matrix{c'_{U^!}:I_o{\buildrel{c_U}\over{\longrightarrow}}U o
U^!{\buildrel{c'_o(U,U^!)}\over{\longrightarrow}}U^! o U}
$$

and

$$
\matrix {d'_{U^!}:U\bullet
U^!{\buildrel{c'_{\bullet}(U,V)}\over{\longrightarrow}}U^!\bullet
U{\buildrel{d_U}\over{\longrightarrow}}I_o}
$$

We have:

$$
\matrix{U^!\rightarrow I_{\bullet}\bullet
U^!{\buildrel{c'_{U^!}\bullet Id_{U^!}}\over{\longrightarrow}}(U^!
o U)\bullet U^!{\buildrel{f_{U^!,U,U^!}}\over{\longrightarrow}}
U^! o (U\bullet U^!){\buildrel{Id_{U^!} o
d'_{U^!}}\over{\longrightarrow}} U^! o I_o}
$$

$$
=\matrix{U^!\longrightarrow I_{\bullet}\bullet
U^!{\buildrel{c_U\bullet Id_{U^!}}\over{\longrightarrow}} (U o
U^!)\bullet U^!{\buildrel{c'_o(U,U^!)\bullet
Id_{U^!}}\over{\longrightarrow}}(U^! o U)\bullet U^!}
$$

$$
\matrix{    {\buildrel{f_{U^!,U,U^!}}\over{\longrightarrow}}U^! o
(U\bullet U^!){\buildrel{Id_{U^!} o
c'_{\bullet}(U,U^!)}\over{\longrightarrow}}U^! o (U^!\bullet
U){\buildrel{Id_{U^!} o d_U}\over{\longrightarrow}} U^! o I_o}
$$

The compatibility property of the quadratic braiding implies:

$$
\matrix{ (U o U^!)\bullet U^!{\buildrel{c'_o(U,U^!)\bullet
Id_{U^!}}\over{\longrightarrow}}(U^! o U)\bullet U^!
{\buildrel{f_{U^!,U,U^!}}\over{\longrightarrow}}U^! o (U\bullet
U^!)}
$$

$$
=\matrix{(U o U^!)\bullet U^!{\buildrel{{c'}_{\bullet}}(U o U^!,
U^!)\over{\longrightarrow}} U^!\bullet (U o
U^!){\buildrel{h_{U^!,U,U^!}}\over{\longrightarrow}}(U^!\bullet U)
o U^!{\buildrel{c'_o(U^!\bullet U,U^!)}\over{\longrightarrow}}U^!
o (U^!\bullet U)}
$$

This implies that:
$$
\matrix{U^!\rightarrow I_{\bullet}\bullet
U^!{\buildrel{c'_{U^!}\bullet Id_{U^!}}\over{\longrightarrow}}(U^!
o U)\bullet U^!{\buildrel{f_{U^!,U,U^!}}\over{\longrightarrow}}
U^! o (U\bullet U^!){\buildrel{Id_{U^!} o
d'_{U^!}}\over{\longrightarrow}} U^! o I_o}
$$

$$
=\matrix{U^!\longrightarrow I_{\bullet}\bullet
U^!{\buildrel{c_U\bullet Id_{U^!}}\over{\longrightarrow}}(U o
U^!)\bullet U^!{\buildrel{{c'}_{\bullet}}(U o U^!,
U^!)\over{\longrightarrow}} U^!\bullet (U o
U^!){\buildrel{h_{U^!,U,U^!}}\over{\longrightarrow}}(U^!\bullet U)
o U^!{\buildrel{c'_o(U^!\bullet U,U^!)}\over{\longrightarrow}}U^!
o (U^!\bullet U)}
$$

$$
\matrix{{\buildrel{Id_{U^!} o
c'_{\bullet}(U,U^!)}\over{\longrightarrow}}U^! o (U^!\bullet
U){\buildrel{Id_{U^!} o d_U}\over{\longrightarrow}} U^! o I_o}
$$

Since $\bullet$ is a tensor product, we deduce that:

$$
\matrix{U^!\longrightarrow I_{\bullet}\bullet
U^!{\buildrel{c_U}\over{\longrightarrow}}(U o U^!)\bullet
U^!{\buildrel{{c'}_{\bullet}(U o
U^!,U^!)}\over{\longrightarrow}}U^!\bullet (U o U^!)}
$$

$$
=\matrix{U^!\longrightarrow U^!\bullet
I_{\bullet}{\buildrel{c_U}\over{\longrightarrow}}U^!\bullet (U o
U^!)}
$$

This implies that

$$
\matrix{U^!\rightarrow I_{\bullet}\bullet
U^!{\buildrel{c'_{U^!}\bullet Id_{U^!}}\over{\longrightarrow}}(U^!
o U)\bullet U^!{\buildrel{f_{U^!,U,U^!}}\over{\longrightarrow}}
U^! o (U\bullet U^!){\buildrel{Id_{U^!} o
d'_{U^!}}\over{\longrightarrow}} U^! o I_o}
$$

$$
=\matrix{U^!\longrightarrow U^!\bullet
I_{\bullet}{\buildrel{Id_{U^!}\bullet
c_U}\over{\longrightarrow}}U^!\bullet (U o
U^!){\buildrel{h_{U^!,U,U^!}}\over{\longrightarrow}}(U^!\bullet U)
o U^!{\buildrel{c'_o(U^!\bullet U,U^!)}\over{\longrightarrow}}U^!
o (U^!\bullet U)}
$$

$$
\matrix{{\buildrel{Id_{U^!} o
c'_{\bullet}(U,U^!)}\over{\longrightarrow}}U^! o (U^!\bullet
U){\buildrel{Id_{U^!} o d_U}\over{\longrightarrow}} U^! o I_o}
$$

Using property $4$ and the fact that $o$ is a tensor product, we
deduce that

$$
\matrix{U^!\rightarrow I_{\bullet}\bullet
U^!{\buildrel{c'_{U^!}\bullet Id_{U^!}}\over{\longrightarrow}}(U^!
o U)\bullet U^!{\buildrel{f_{U^!,U,U^!}}\over{\longrightarrow}}
U^! o (U\bullet U^!){\buildrel{Id_{U^!} o
d'_{U^!}}\over{\longrightarrow}} U^! o I_o}
$$

is the identity of $U^!$.
\bigskip

We also have:

$$
\matrix{U\longrightarrow U\bullet
I_{\bullet}{\buildrel{Id_U\bullet
c'_{U^!}}\over{\longrightarrow}}U\bullet(U^! o
U){\buildrel{h_{U,U^!,U}}\over{\longrightarrow}} (U\bullet U^!) o
U{\buildrel{d'_{U^!} o Id_U}\over{\longrightarrow}}U o
I_o\rightarrow U}
$$

$$
=\matrix{U\longrightarrow U\bullet
I_{\bullet}{\buildrel{Id_U\bullet c_U}\over{\longrightarrow}}
U\bullet(U o U^!){\buildrel{Id_U\bullet c'_o(U,U^!)}\over
{\longrightarrow}} U \bullet (U^! o U)}
$$

$$
\matrix{{\buildrel{h_{U,U^!,U}}\over{\longrightarrow}}(U\bullet
U^!) o U{\buildrel{c'_{\bullet}(U,U^!) o
Id_U}\over{\longrightarrow}}(U^!\bullet U) o U{\buildrel{d_U o
Id_U}\over{\longrightarrow}} U o I_o}
$$

Using the compatibility property of the braiding, we have:

$h_{U,U^!,U}:U \bullet (U^! o U)\longrightarrow (U \bullet U^!) o
U $

$$
=\matrix{U \bullet (U^! o U){\buildrel{c'_{\bullet}(U, U^! o
U)}\over{\longrightarrow}} (U^! o U)\bullet
U{\buildrel{c'_o(U^!,U)\bullet Id_U}\over{\longrightarrow}}(U o
U^!)\bullet U{\buildrel{f_{U,U^!,U}}\over{\longrightarrow}}U o
(U^!\bullet U){\buildrel{c'_o(U,U^!\bullet
U)}\over{\longrightarrow}}(U^!\bullet U) o U}
$$

This implies that:

$$
\matrix{U\longrightarrow U\bullet
I_{\bullet}{\buildrel{Id_U\bullet
c'_{U^!}}\over{\longrightarrow}}U\bullet(U^! o
U){\buildrel{h_{U,U^!,U}}\over{\longrightarrow}} (U\bullet U^!) o
U{\buildrel{d'_{U^!} o Id_U}\over{\longrightarrow}}U o
I_o\rightarrow U}
$$

$$
=\matrix{U\longrightarrow U\bullet
I_{\bullet}{\buildrel{Id_U\bullet c_U}\over{\longrightarrow}}
U\bullet(U o U^!){\buildrel{Id_U\bullet c'_o(U,U^!)}\over
{\longrightarrow}} U \bullet (U^! o U)}
$$

$$
\matrix{{\buildrel{c'_{\bullet}(U, U^! o
U)}\over{\longrightarrow}} (U^! o U)\bullet
U{\buildrel{c'_o(U^!,U)\bullet Id_U}\over{\longrightarrow}}(U o
U^!)\bullet U{\buildrel{f_{U,U^!,U}}\over{\longrightarrow}}U o
(U^!\bullet U){\buildrel{c'_o(U,U^!\bullet
U)}\over{\longrightarrow}}(U^!\bullet U) o U{\buildrel{d_U o
Id_U}\over{\longrightarrow}} U o I_o}
$$

Using the fact that $\bullet$ is a tensor product, we deduce that:

$$
\matrix{U\longrightarrow U\bullet
I_{\bullet}{\buildrel{Id_U\bullet c_U}\over{\longrightarrow}}
U\bullet(U o U^!){\buildrel{c'_{\bullet}(U,U o U^!)}\over
{\longrightarrow}} (U o U^!)\bullet U}
$$

$$
=\matrix{U\longrightarrow I_{\bullet}\bullet
U{\buildrel{c_U}\over{\longrightarrow}}(U o U^!)\bullet U}
$$

This implies that:

$$
\matrix{U\longrightarrow U\bullet
I_{\bullet}{\buildrel{Id_U\bullet
c'_{U^!}}\over{\longrightarrow}}U\bullet(U^! o
U){\buildrel{h_{U,U^!,U}}\over{\longrightarrow}} (U\bullet U^!) o
U{\buildrel{d'_{U^!} o Id_U}\over{\longrightarrow}}U o
I_o\rightarrow U}
$$

$$
=\matrix{U\longrightarrow I_{\bullet}\bullet
U{\buildrel{c_U}\over{\longrightarrow}}(U o U^!)\bullet U
{\buildrel{f_{U,U^!,U}}\over{\longrightarrow}}U o (U^!\bullet
U){\buildrel{c'_o(U,U^!\bullet
U)}\over{\longrightarrow}}(U^!\bullet U) o U{\buildrel{d_U o
Id_U}\over{\longrightarrow}} U o I_o}
$$

Using property $3$, and the fact that $o$ is a tensor product we
deduce that

$$
\matrix{U\longrightarrow U\bullet
I_{\bullet}{\buildrel{Id_U\bullet
c'_{U^!}}\over{\longrightarrow}}U\bullet(U^! o
U){\buildrel{h_{U,U^!,U}}\over{\longrightarrow}} (U\bullet U^!) o
U{\buildrel{d'_{U^!} o Id_U}\over{\longrightarrow}}U o
I_o\rightarrow U}
$$

is the identity.

Since the quadratic dual is unique up to an isomorphism, we deduce
that $U$ is a quadratic dual of $U^!$.

\bigskip

{\bf Proposition 9.}

{\it Let $(C,o,\bullet,!)$ be a quadratic category endowed with a
quadratic braiding. The opposite functor $P:U\rightarrow U^!$ is
an equivalence of category.}

\medskip

{\bf Proof.}

Let $U_1$, and $U_2$ be objects of $C$. Using the previous result,
we can suppose that ${U^!_1}^!$ and ${U_2^!}^!=U_2$. Using the
theorem $2$, we deduce a bijection between
$Hom_C(U_1\bullet{U_2^!},I_o)$ and $Hom_C(U_1,U_2)$, and an
isomorphism between $Hom_C(U^!_2\bullet U_1)$ and
$Hom_C(U^!_2,U^!_1)$. The commutative constraint $c'_{\bullet}$
defines an isomorphism between $Hom_C(U_1\bullet U^!_2,I_o)$ and
$Hom_C(U^!_2\bullet U_1,I_o)$. The object $U_1^!$ represents the
functor $n_{U_1}:U\rightarrow Hom_C(V\bullet U_1,I_o)$. This
implies that the correspondence defined on $C$ by $U\rightarrow
U^!$ is functorial, and using the Yoneda lemma, we deduce that the
morphisms between $n_{U_2}$ and $n_{U_1}$ are given by
$Hom_C(U_2^!\bullet U_1,I_o)$. This implies that $P$ is fully
faithful$\bullet$

\bigskip

{\bf Definition 10.}

Let $(C,o,\bullet,!)$ be a quadratic category endowed with a
quadratic braiding. We denote by $C'$ the subcategory of $C$, such
that for each object $U$ of $C'$, there exists an object $V$ of
$C$ such that $U$ is isomorphic to $V o I_{\bullet}$, $(C',o)$ is
a subtensor category of $(C,\bullet)$. We say that the quadratic
category $(C,o,\bullet,I_o)$ is quadratic rigid if for every
objects, $U_1$ and $U_2$ of $C'$, $U_1 o U_2=U_1\bullet U_2$ and
the restriction of the braided associative constraints
$f_{U_1,U_2,U_3}$ and $h_{U_1,U_2,U_3}$ to $C'$ coincide with the
associative constraint of $o$ and $\bullet$.

\medskip

Let $(C,o,\bullet,!)$ be a rigid quadratic category, and  $U$ an
object of $C'$, we define $U^*$ to be $U^! o I_{\bullet}$. We have
$I_{\bullet}=I_{\bullet} o I_o$. This implies that $I_{\bullet}$
is an object of $C'$. Thus $I_{\bullet}\bullet
I_{\bullet}=I_{\bullet} o I_{\bullet}=I_{\bullet}$. We deduce that
for each object $U$ of $C'$, $U o U^!$ is isomorphic to $(U o
I_{\bullet}) o U^!$ which is isomorphic to $U o (U^! o
I_{\bullet})$ Since $U$ is  an object of $C'$, $U o U^!$ is
isomorphic to $U\bullet (U^! o I_{\bullet})=U\bullet U^*$.

\bigskip

{\bf Proposition 11.}

{\it  Let $U$ be an object of $C'$ the map $c_U$, and the map

$$
\matrix{d^1_{U^!}:U o(U^! o I_{\bullet}) =U\bullet (U^! o
I_{\bullet}){\buildrel{h_{U,U^!,I_{\bullet}}}\over{\longrightarrow}}(U
\bullet U^!) o I_{\bullet}{\buildrel{d'_U o
I_{\bullet}}\over{\longrightarrow}} I_o o I_{\bullet}=I_{\bullet}}
$$

defines on $(C',\bullet)$ the structure of a rigid tensor
category.}

\medskip
{\bf Proof.}

We have to show  that for each element $U$ in $C'$:

$$
\matrix{U{\buildrel{c'_U}\over{\longrightarrow}}U\bullet(U^! o
U)=U\bullet (U^*\bullet
U){\buildrel{h_{U,U^!,U}}\over{\longrightarrow}}(U\bullet U^*) o U
{\buildrel{Id_U\bullet d^1_U }\over{\longrightarrow}}U}
$$
is the identity of $U$, and

$$
\matrix{U^*{\buildrel{c'_U}\over{\longrightarrow}}(U^! o U)\bullet
U^*=(U^*\bullet U)\bullet
U^*{\buildrel{f_{U^!,U,U^*}}\over{\longrightarrow}}U^* o(U\bullet
U^*){\buildrel{Id_{U^*} o d^1_U}\over{\longrightarrow}}U^*}
$$
is the identity. These assertions follows from the fact that $U$
is the dual of $U^!$. In the second assertion, we multiply $(4)$
applying to $U^!$ by $I_{\bullet}$. This implies that $U^*$ is a
dual of $U$ since the category is braided we deduce from Deligne
that this category is rigid$\bullet$

\bigskip

{\bf Definition 12.}

Let $(C,o,\bullet, !)$ be a quadratic rigid category, and
$U\rightarrow U$ a morphism. We can define $h o Id_{I_{\bullet}}:U
o I_{\bullet}\rightarrow U o I_{\bullet}$. We define $Trace(h)$ to
be the trace of $h o Id_{I_{\bullet}}$. This is the endomorphism
of $I_{\bullet}$ defined by:

$$
\matrix{I_{\bullet}{\buildrel{c'_U}\over{\longrightarrow}} U^! o
U{\buildrel{Id_U^! o h}\over{\longrightarrow}}U^! o
U{\buildrel{d^1_U}\over{\longrightarrow}}I_{\bullet}}
$$

We denote by $rank(U)$ the trace of $Id_U$. Let $h,h':U\rightarrow
U$, $Trace(hh')= Trace(h)Trace(h')$.

The ring $Hom(I_{\bullet},I_{\bullet})$ is commutative since
$I_{\bullet}$ is the neutral element of $(C,\bullet)$.

\bigskip

{\bf Proposition 13.}

{\it Let $(C,o,\bullet,!)$ be a quadratic rigid category, and $C"$
be the subcategory of $C$ such that for every object $U$ of $C"$,
there exists an object $V$ of $C$ such that $U=V\bullet I_o$. The
category $(C",o,I_o)$ is a rigid tensor category.}

\medskip

{\bf Proof.}

The restriction of the opposite functor $P:C\rightarrow C$ defined
on $C$ by $U\rightarrow U^!$ to $C'$ define an isomorphism of
tensor categories between $C'$ and $C"$.

\bigskip

{\bf Definition 14.}

 Let $(C,o,\bullet,!)$ a quadratic category, and $h:U\rightarrow
V$ a morphism. The morphism $h':U^!\rightarrow V^!$ is a
contragredient of $h$ if and only if $h o h'\circ c_U=c_V$, and
$d_V\circ (h'\bullet h)=d_U$.

\bigskip

{\bf Proposition 15.}

{\it Let $(C,o,\bullet,!)$ be a braided quadratic category. Then a
contragredient map $h$ is invertible.}

\medskip

{\bf Proof.}

Let $h:U\rightarrow V$ be a contragredient map defined by the
application $h':U^!\rightarrow V^!$. We are going to show that
$h'\circ h^!=Id_{V}$  since the category is braided, this implies
that $h$ is invertible.

$h^!\circ h'=$
$$
\matrix{U^!{\buildrel{h'}\over{\longrightarrow}}V^!{\buildrel{Id_{V^!}o
c_U}\over{\longrightarrow}}V^!\bullet(U o
U^!){\buildrel{h_{U^!,U,U^!}}\over{\longrightarrow}}(V^!\bullet U)
o U^!{\buildrel{(Id_{V^!}\bullet h) o
Id_{U}}\over{\longrightarrow}}(V^!\bullet V) o
U^!{\buildrel{d_V}\over{\longrightarrow}}U^!}
$$

Using property $5$, we deduce that $h^!\circ h'=$

$$
\matrix{U^!{\buildrel{c_U}\over{\longrightarrow}}U^!\bullet(U o
U^!){\buildrel{h_{U^!,U,U^!}}\over{\longrightarrow}}(U^!\bullet U)
o U{\buildrel{h'\bullet h}\over{\longrightarrow}}(V^!\bullet V)  o
U^!{\buildrel{d_V}\over{\longrightarrow}} U^!}
$$

Using the fact that $d_V\circ (h'\bullet h)=d_U$, and property
$4$, we deduce that $h^!\circ h'=Id_{U^!}$ $\bullet$

\bigskip

 {\bf Definition 16.}

A quadratic functor $H:(C,o,\bullet,!)\rightarrow
(C',o',\bullet',!')$ is a functor  of tensor categories
$H:(C,o)\rightarrow (C',o')$ which commutes with the quadratic
dual, $H(U^!)=H(U)^!.$

\bigskip

{\bf Proposition 17.}

{\it Let $(C,o,\bullet,!)$ and $(C',o',\bullet',!')$ be quadratic
baided categories, and $F,F':(C,o,\bullet,!)\longrightarrow
(C',o',\bullet',!')$ be two quadratic functors every morphism
$u:F\rightarrow F'$ is an isomorphism.}

{\bf proof.}

The morphism $u:F\rightarrow F'$ is defined by a family of
morphisms $u_U:F(U)\rightarrow F'(U)$ where $U$ is an object of
$C$. The morphism $u_{U^!}$ is a contragredient of $u_U$. This
implies that $u$ is an isomorphism. The inverse of $u$ is the
family of maps $u^{!'}_{U^!}$.

\bigskip

We have the morphism $c_{I_o}:I_o\rightarrow I_o o
I^!_o=I_{\bullet}$.

\bigskip

Let $U$ and $V$ be objects of $C$, we have a morphism $U\bullet
V\longrightarrow U o V$ defined by

$$
\matrix{U\bullet (V o I_o){\buildrel{c'_{\bullet}(U, V o
I_o)}\over{\longrightarrow}}(V o I_o)\bullet
U{\buildrel{f_{V,I_o,U}}\over{\longrightarrow}}V o (I_o\bullet
U){\buildrel{Id_V o(c_{I_o}\bullet Id_U)}\over{\longrightarrow}}V
o (I_{\bullet}\bullet
U){\buildrel{c'_o(V,U)}\over{\longrightarrow}}U o V}
$$

\medskip

\bigskip

{\bf II. Koszul complexes.}

\bigskip

\medskip

We suppose that the category $C$ is additive, and is contained in
an abelian category $C'$.
  Let $h$ be an element of $Hom(U,U)$, the
isomorphism $Hom(U,U)\rightarrow Hom(I_{\bullet},U o U^!)$,
defines a map $d'_h:I_{\bullet}\rightarrow U o U^!$. We can
construct the map:

$$
\matrix{d_h: U o U^!\longrightarrow I_{\bullet}\bullet (U o
U^!){\buildrel{d'_h\bullet Id_{U o U^!}}\over{\longrightarrow}}(U
o U^!)\bullet (U o U^!){\buildrel{l_U}\over{\longrightarrow}}U o
U^!}
$$

 \bigskip

 {\bf Definitions 1.}

 The category $C$ is an $n$-Koszul category, if for every object
 $U$, and each map $h:U\rightarrow U$, $(d_h)^{n+1}=0$. The
 category $C$ is a Koszul category if it is $1$-Koszul. This is equivalent
 to saying that  $(d_h)^2=0$.

\bigskip
Let $C$ be a Koszul category. We denote by $D_U$ the
endomorphism $d_{Id_U}$, and define the first Koszul complex to be
$L(U)=(U o U^!,D_U)$. We say that  $U$ is Koszul if $L(U)$ is
exact. Let

$$
\matrix{U{\buildrel{d_0}\over{\longrightarrow}}U_1...U_p
{\buildrel{d_p}\over{\longrightarrow}}U_{p+1}...}
$$

be a resolution of $U$ in $C'$. We say that this resolution is a
Koszul resolution, if there exists an object $V$ of $C$ endowed
with a differential $\alpha_V$, such that there exists  embedding
$\alpha_U: U o U^!\longrightarrow V$, $e_p:U_p\longrightarrow V$,
such that the following squares are commutative:

$$
\matrix{U_p{\buildrel{d_p}\over{\longrightarrow}}U_{p+1}\cr
\downarrow e_p\downarrow e_{p+1}\cr
V{\buildrel{\alpha_V}\over{\longrightarrow}}V}
$$

$$
\matrix{U o U^!{\buildrel{\alpha_U}\over{\longrightarrow}}V\cr
\downarrow D_U\downarrow \alpha_V\cr U o
U^!{\buildrel{\alpha_U}\over{\longrightarrow}}V}
$$
\bigskip

{\bf  Koszul complexes of algebras.}

\bigskip

We suppose that $C$ is a Koszul category, the objects of $C$ are
graded algebras defined over a field $F$. We denote by $U_i$ the
$i$-component of $U$, and we suppose that $U_0=F$.   The map $D_U$
is a left multiplication by an element $\alpha'_U$ of $U o U^!$,
we suppose that $C$ is stable by the usual tensor product of
$F$-vector spaces,and there exists an embedding $U o V\rightarrow
U\otimes V$. We denote by $\alpha_U$ the image of $\alpha'_U$ by
this embedding. We suppose that $\alpha_U\in U_1\otimes U^!_1$.
The family $(U\otimes U^!_i,\alpha_U)$ is a Koszul complex called
the first algebra Koszul complex. If this complex is exact, the
algebra is called a Koszul algebra.

\bigskip

{\bf Theorem 2.}

{\it Let $U$ be an object of $C$, if
 the complex $(U\otimes U^!,\alpha_U)$ is acyclic then
 $Ext_{U}(F,F)=U^!$.}

\bigskip

{\bf Proof.}

 Suppose that $(U \otimes U^!,\alpha_U)$ is acyclic. This implies that:

$$
0\longrightarrow U\longrightarrow U\otimes U^!_1...\longrightarrow
U\otimes {U^!}_i\longrightarrow U\otimes{U^!}_{i+1}...
$$

is a $U$-resolution of $U$. The tensor product $(U\otimes_F
{U^!}_i)\otimes_U F$ is isomorphic to $(U\otimes_U
F)\otimes_F{U^!}^i$. The tensor product $U\otimes_U F$ is $F$,
since the left-module of $F$ verifies $I(U)F=0$ where $I(U)$ is
the augmentation ideal. This implies that$(U\otimes_U
F)\otimes_F{U^!}^i={U^!}_i$ and the multiplication by $\alpha_U$
induces the zero map between $U^!_i$ and $U^!_{i+1}$. This implies
that $Ext_U(F,F)=U^!$$\bullet$

\bigskip

{\bf Second Koszul complex.}

\medskip

Let ${U^!}^*$ be the $F$-dual of $U^!$, $U\otimes {U^!}^*$ is
embedded in $Hom_U(U\otimes U^!,U)$ as follows: let $u_1,u_2$ be
elements of $U$, $v_1$ an element of ${U_l^!}^*$ and $v_2$ and
element of $U^!_p$. We define $(u_1\otimes v_1)(u_2\otimes
v_2)=v_1(v_2)u_1u_2$ if $l=p$, $(u_1\otimes v_1)(u_2\otimes
v_2)=0$ if $l\neq p$. We can define the differential $D'_U$ on
$U\otimes {U^!}^*$ by setting $D'_U(h)(u)=h(\alpha_U(u))$. The
complex $(U\otimes U^!,\alpha_U)$ and $(U\otimes {U^!}^*,D'_U)$
are dual each other.

\bigskip

{\bf Quadratic algebras.}

\bigskip

A quadratic algebra $U$ is an homogeneous algebra generated by its
elements of degree $1$, equivalently, there exists a vector space
$V$, such that $U$ is the quotient of the tensor product $T(V)$ by
an homogeneous ideal $C$ generated by a subspace contained in
$V\otimes V$. The category of quadratic algebras is a quadratic
category, Let $V^*$ be the dual space of $V$ one defines: $U^!$ to
be the quotient of $T(V^*)$ by $C^{\bot}$ where $C^{\bot}$ is the
annihilator of $C$.

Suppose that $U'=T(V')/C,$ is another quadratic algebra, we denote
by $t_{23}:V^{\otimes^2}\otimes {V'}^{\otimes^2}\rightarrow
V\otimes V'\otimes V\otimes V'$ the map defined by
$t_{23}(u_1\otimes u_2\otimes v_1\otimes v_2)=u_1\otimes
v_1\otimes u_2\otimes v_2$. Then one defines

$U\bullet U'=T(V\otimes V')/t_{23}(C\otimes C')$,

$U o U'= T(V\otimes V')/t_{23}(V^{\otimes^2}\otimes C'\oplus
C\otimes {V'}^{\otimes^2})$.

One of the main objective of specialists in Koszul structures is
to construct a Koszul resolution of an object $U$. The reason of
this is the fact that using this resolution we  can  easily
compute $U$-homology and cohomology.
 The usual Bar complex allows to compute the homology of an algebra.
 Recall its definition. Let $U$ be an algebra, $\epsilon:U\rightarrow F$ the
augmentation of $U$, and $I(U)$ the kernel of $\epsilon$. The Bar
complex of $U$ is the tensor product $B(U)=U\otimes T(I(U))\otimes
U$. Its elements are denoted by $u[u_1:...:u_p]u'$. The Bar
complex is bigraded the degree of $v=u[u_1:...:u_p]u'$ is $(m,n)$
where $m=\sum_{i=1}^{i=p}degree(u_i)+degree(u)+degree(u')$, and
$n=p$ is the homological degree.
 The converse of the theorem above  for quadratic algebras is shown by Priddy using a
 spectral sequence. Here is an elementary proof:

 \bigskip

 {\bf Proposition 3.}

 {\it Let $U$ be a quadratic algebra, if $Ext_U(F,F)=U^!$, then
 $U$ is a Koszul algebra.}

\medskip

{\bf Proof.}

Suppose that $U=T(V)/C$. We compute $Ext_U(F,F)$ using the bar
complex. We have $F\otimes B(U)\otimes F=T(I(U))$. This implies
that $Ext_U(F,F)=H^*(T(I(U))^*)$, since $U\otimes_UF=F$. We denote
by $D'$ the differential of this complex. The algebra $T(V^*)$ is
contained in $T(I(U)^*)$ since $I(U)$ contains $V$. We use the
homological degree to define the graduation of $T(I(U)^*)$. We
obtain that $T(I(U)^*)_0=T(V^*)$, and $T(I(U))_1=\sum
(V^{\otimes^p}\otimes (T(V)/C)\otimes V^{\otimes^l})^*$. The
kernel of the restriction of $D'$ to $T(V^*)$ is $U^!$. The
complex $(I(U)\otimes {U^!}^*,D'_U)$ is contained in the Bar
complex $T((I(U),D)$. The result follows from the fact that
$H^{p,l}(T(I(U))=H_{p,l}(T(I(U))=0$ if $l\neq p$.

\bigskip
\bigskip

{\bf III. Quadratic category of operads and $n$-categories.}

\bigskip

\medskip

The theory of operads represents the mathematical framework to
encodes operations. The binary operations are encoded by by
quotient of free operads of rank $2$. These operads are defined
and studied by Ginzburg and Kapranov. An important problem still
unsolved in mathematics is to define a theory of $n$-categories.
On this purpose one needs to define the coherence relations for
$n$-arrows. In this part, we will define the notion of
$n$-quadratic operad, and show how a $n$-coherence system of a
theory of $n$-category is related to this notion. The notion of
$n$-Koszul algebras has been recently defined by Berger, [3]..

\bigskip

 {\bf Definition 1.}

Let $S_n$ be the symmetric group of $n$ elements, $d\in S_n$,
$d'\in S_{n'}$ two permutations. For every $p<n+1$, we define the
permutation $d o_p d'\in S_{n+n'-1}$ by:

Suppose that $j< p$ or $j>p+n$, $(d c_p d')(j)= d(j)$ if $d(j)<
d(p)$ otherwise $(d c_p d')(j)=d(j)+n'-1$.

Suppose that $p-1<j<p+n'-1$, then $(d c_p d')(j)=d(j)+d'(j-p+1)$.

A linear operad is a family of vectors spaces $(P(n))_{n\in{\N}}$,
such that for each  $p<n+1$, there exists a composition
$c_p:P(n)\otimes P(n')\longrightarrow P(n+n'-1)$ which satisfies
the following conditions:

$P(n)$ is endowed with an action of $S_n$, such that for every
elements $u_n\in P(n)$, $u'_{n'}\in P(n')$, $d\in S_{n'}$, and
$d'\in S_{n'}$,

$$
(u_n c_p u'_{n'})^{d c_p d'}={u_n}^d c_{d(p)} {u'_{n'}}^{d'}.
$$

Let $u_1\in P(n_1)$, $u_2\in P(n_2)$, and $u_3\in P(n_3)$, given
integer $0<p<p'<n_1+1$,

$$
(u_1 c_p u_2) c_{p'+n_2-1} u_3=(u_1 c_{p'} u_3) c_{p} u_2
$$

Suppose that $0<p<n_1+1$, and $0<p'<n_2+1$

$$
(u_1 c_p u_2) c_{p+p'-1} u_3 = u_1 c_{p} (u_2 c_{p'} u_3)
$$

\bigskip

We define the fundamental example of operad which allow to define
the notion of algebra.

\bigskip

{\bf Definition 2.}

Let $V$ be a vector space, we denote by
$P(n)=Hom(V^{\otimes^n},V)$ the space of linear maps
$V^{\otimes^n}\longrightarrow V$. The family of vector spaces
$(P(n))_{n\in {\N}}$ defines an operad $F(V)$ such that for every
element $h\in P(n)$, and $h'\in P(n')$ and $p\leq n$,

$$
(h c_p h')(v_1\otimes..\otimes
v_{n+n'-1})=h(v_1\otimes..h'(v_p\otimes..\otimes v_{p+n'-1})
..\otimes v_{n+n'-1}).
$$

\bigskip

We present now the notion of free operad which will allow us to
define the notion of $n$-Koszul operad. Let $P=(P(n))_{n\in {\N}}$
be a family of vector spaces such that each $P_n$ is endowed with
an action of $S_n$. Consider the category $C_P$, whose objects are
morphisms $h:P\longrightarrow P'$, where $P'$ is an operad, $h$ is
a collection of morphisms $h_n:P(n)\longrightarrow P'(n)$ which
commute with the action of $S_n$. The category $C_P$ has an
initial object called the free operad generated by $P$. Recall the
construction of this operad.

\bigskip

{\bf Definition 3.}

Consider a tree $T$,  $V$ the set of its vertex. Let $v$ be an
element of $V$, we denote by $In(v)$ the set of inner edges of
$v$. If the cardinal of $In(v)$ is $p$, we define $P(In(v))=E(p)$.
We defined $E(T)$ to be the tensor product $\otimes_{v\in
V}P(In(v))$, and by $F(n)=\oplus_{n-tree}E(T)$. The family of
vector spaces $F(n)$ is the free operad generated by $(E(n))_{n\in
{\N}}$. The composition $c_p$ are induced by the corresponding
composition of trees.

An ideal $U$ of an operad $P$, is a family of subspaces $U(n)$
included in $P(n)$, such that for every elements $u_{n_1}\in
P(n_1)$, $u_{n_2}\in P(n_2)$, $v_{n_1 }\in U(n'_1)$, $v'_{n'_2}\in
U(n'_2)$, $p\leq n_1$, $p'\leq n'_2$, $u_{n_1} c_p v_{n_1}\in
U(n_1+n'_1-1)$, and $v_{n_2} c_{p'} u_{n_2}\in U(n_2+n'_2-1)$. We
can construct the quotient of $P$ by $U$.

An $n$-quadratic operad, is an operad $P$, which is the quotient
of a free operad $F(E)$ such that $E(n)=0$ if $n\neq 2$ by an
ideal generated by a subspace of $F(E)(n)$. A morphism
$h:U\rightarrow U'$ between two $n$-quadratic operads $F(E)/C$ and
$F(E')/C'$ is a morphism of operad generated by a morphism of
$S_2$-spaces $h_2:E(2)\longrightarrow E'(2)$.

\bigskip

{\bf Examples.}

\bigskip

{\bf The operad $As$.}

\bigskip

The operad $As$ is the operad which encodes the associativity.
This operad is defined as follows: Consider the family of vectors
spaces $(E(n))_{n\in{\N}}$ such that $E(n)=0$ if $n\neq 0$, $E(2)$
is the $2$-dimensional vector space whose generators are denoted
by $u_1u_2$, $u_2u_1$. We define on $E(2)$ the action of $S_2$
defined by $h(u_1u_2)=u_2u_1$. The free operad $F(E)$ generated by
$(E(n))_{n\in{\N}}$ verifies $F(E)(1)=0$, $F(E)(2)$ is $3$-copies
of $E(2)\otimes E(2)$ denoted by $3E(2)\otimes E(2)$ corresponding
to the trees with one inner arrow. We will denote $(u_1u_2) c_1
(u_1u_2)$ by $(u_1u_2)u_3$, $(u_1u_2) c_1 (u_2u_1)$ by
$u_3(u_1u_2)$, $(u_1u_2) c_2 (u_1u_2)=u_1(u_2u_3)$,... The
operation $(u_1u_2) c_j (u_1u_2)$ is obtained by the putting a
tree with two edges on the $j$-edge of a tree with $2$-edges,
$u_1$, $u_2$, $u_3$ are put on the outer vertices from left to
right. We define the ideal $J_{As}$ to be the ideal generated
$u_{i_1}(u_{i_2}u_{i_3})-(u_{i_1}u_{i_2})u_{i_3}$, where
$u_{i_1}(u_{i_2}u_{i_3})$ and $(u_{i_1}u_{i_2})u_{i_3}$ are chosen
in the generators of $3E(2)\otimes E(2)$ described above. The
quotient of $F(E)$ by $J_{As}$ is the associative operad.

An associative structure on the vector space $V$ is defined by a
morphism of operads between $As\rightarrow F(V)$.

\bigskip

{\bf Coherence relations in $n$-category and $n$-quadratic
operads.}

\bigskip

The theory of $n$-categories is not well-understood, to define
such a theory, one need to describe the coherence relations
satisfied by $n$-arrows. The operad $As$ describe the coherence
relation for arrows in the theory of category. Otherwise said the
composition of arrows in this theory is  associative. We construct
a $2$-quadratic operad which describe the coherence of $2$-arrows
in the theory of $2$-categories called the pentagon operad and
denoted $P_2$. A theory of $n$-categories using operads has been
studied  by Batanin.

\bigskip

{\bf The Pentagon operad.}

\bigskip

To define the Pentagon operad, we consider the free operad
generated by the family of vector spaces $(E(n))_{n\in {\N}}$ such
that $E_n=0$ if $n\neq 2$, and $E(2)$ is the linear space
generated by $u_1u_2, u_2u_1$. The generators of $F(E)(4)$
$((u_1u_2)u_3)) c_1 (u_1u_2)$ is denoted by $((u_1u_2)u_3)u_4$,
$((u_1u_2)u_3) c_2 (u_1u_2)$ by $u_1((u_2u_3)u_4)$, $(u_1u_2) c_3
((u_1u_2)u_3)=(u_1u_2)(u_3u_4)$,...

The ideal $J_{P_2}$ defining the pentagon operad is generated by
$(u_{i_1}(u_{i_2}u_{i_3}))u_{i_4}-((u_{i_1}u_{i_2})u_{i_3})u_{i_4}+
(u_{i_1}u_{i_2})(u_{i_3}u_{i_4})+u_{i_1}(u_{i_2}(u_{i_3}u_{i_4}))-
u_{i_1}((u_{i_2}u_{i_3})u_{i_4})$. The operad $P_2$ is the
quotient of $F(E)$ by $P_2$.

\bigskip

{\bf The quadratic operad $P_n$.}

\bigskip

Supposed defined the operad $P_n$ which describe the coherence
relations of the $n$-arrows of an $n$-category. It is the quotient
of the free operad $F(E)$ considered to defined the operad $As$,
by an ideal $J_{P_n}$ generated by elements of $F(E)(n+2)$. The
generators of $J_{P_n}$ are elements $u_n-u'_n$ where $u_n$
represents a morphism $h_n$ between
$u_{i_1}(u_{i_2}(...(u_{n+1}u_{i_{n+2}}))$ and
$((u_{i_1}u_{i_2})u_{i_3})...)u_{n+2})$.

We can define the generators of $F(E)(n+3)$ recursively as we did
for the $As$ and $P_2$ operad, their are of the form
$u_1...u_{n+3}$ we do not precise what are the brackets. A path
$(v_1=u_{i_1}(...(u_{i_{n+2}}u_{i_{n+3}})),u_2,...,v_p=
((u_{i_1}u_{i_2})...)u_{n+2})u_{n+3})$ is a family of elements
$v_j$ between the set of generators described above such that
there exists a morphism between $h_j:v_j\rightarrow v_{j+1}$ which
occurs in the definition of $J_{P_l}, l\leq n$ (we have to make
blocks to identified the $v_j$ to elements of $F(E)(n+2)$ for
example $u_1(u_2u_3)$ is identified to $u_1v'_1$ where
$v'_1=u_2u_3$). The ideal $J_{P_{n+1}}$ is the ideal generated by
the sets $v_1+...+v_p-(v'_1+...+v'_{p'})$ where $(v_1,..,v_p)$ and
$(v'_1,..,v'_n)$ are paths. We define the operad $P_{n+1}$ to be
the quotient of $F(E)$ by $J_{P_{n+1}}$.

\bigskip

{\bf Dual of a $n$-quadratic operad.}

\bigskip

{\bf Definition 5.}

Let $E(n), n\in {\N}$ be a family of vector spaces such that
$E(n)$ is endowed with an action of $S_n$, we define the
$S_n$-module $E(n)^{\vee}$ as follows, for each element $v$ in the
dual ${E(n)}^{*}$ of $E(n)$, $v'\in E(n)$, and $u\in S_n$,
$(u(v))(v')=sign(u)v(u(v'))$.

Let $U=F(E)/C$ be a $n$-quadratic operad. We define the dual $U^!$
to be the  quotient of the free operad generated by the family
$({E(p)}^{\vee})_{p\in{\N}}$ by the annihilator $C^{\bot}$ of $C$.

\bigskip

 The operad ${U^!}^!$ is $U$. Let $h:U\longrightarrow U'$ be a
 morphism of quadratic operads generated by the $S_2$-map
 $h_2:E(2)\rightarrow E'(2)$. The dual of $h_2$ defines a morphism
 $h^!{:U'}^!\rightarrow U^!$, and the correspondence
 $U\longrightarrow U^!$ is an opposite functor.

\bigskip

{\bf Tensor products in the category of operads.}

\bigskip

{\bf Definition 6.}

 Let $E(2)$ and $E'(2)$ be two vector spaces, $F(E)$ and $F(E')$ the free
 operads respectively generated by $E(2)$ and $E'(2)$ and
  $P=F(E)/C$, and $P'=F(E')/C'$ be two $n$-quadratic operads,
  consider the free operad generated by $E(2)\otimes E(2)'$.
We can write The usual tensor product $P\otimes P'$ is an operad.
We denote by $P o P'$. The suboperad of $P\otimes P'$ generated by
$E(2)\otimes E(2)'$, and by $P\bullet P'$ the operad $(P^! o
{P'}^!)^!$.

\bigskip

{\bf Proposition 1.}

{\it Let $I_o$ be  the free binary operad generated by the ground
field $L$. Then $I_o$ is the neutral element of $o$,
$I^!_o=I_{\bullet}$ is the neutral element of $\bullet$.}

\bigskip

{\bf Proof.}

Let $P$ be a $n$-quadratic operad quotient of $F(E)$ by the ideal
$C$. The tensor product $E(2)\otimes L$ is $E(2)$. We have $(P o
I_o)(n) =\oplus_{T n tree}\otimes_{vertice T}u\otimes v$ where $u$
is an element of $P$ and $v$ an element of $L$. This implies that
$I_o$ is the neutral element of $P$.  The fact that $I_{\bullet}$
is the neutral element of $\bullet$ follows by duality.

\bigskip

Let $P=F(E)/C$ be a quadratic operad. The evaluation
$E(2)^{\vee}\otimes E(2)\longrightarrow L$ defines a map $d_U: U^!
\bullet U\longrightarrow I_o$ since the image of the ideal which
defines $U^!\bullet U$ is by definition zero. We deduce the
existence of a map $c_U:I_{\bullet}\longrightarrow U o U^!$. The
category $(C, o,\bullet,!)$ whose objects are $n$-quadratic
operads is a quadratic category. We deduce that for every
$n$-quadratic operads $U_1$ and $U_2$, the functor
$U\longrightarrow Hom(U_1\bullet U,U_2)$ is corepresentable by
$\underline{Hom}(U_1,U_2)=U_2 o U^!_1$. And there exists a product
$l_U$ on $U o U^!$ defines by:

$$
\matrix{l_U: U o U^!\longrightarrow I_{\bullet}\bullet (U o
U^!){\buildrel{c_U\bullet Id_{U o U^!}}\over{\longrightarrow}}(U o
U^!)\bullet (U o U^!){\buildrel{l_U}\over{\longrightarrow}}U o
U^!}
$$

The morphism $c_U$ is defined by the map $L\longrightarrow
E(2)\otimes E(2)^{\vee}$ $1\longrightarrow
\delta_U=\sum_{i=1}^{i=n}u_i\otimes u_i^*$, where
$(u_i)_{i=1,..,p}$ is a basis of $E(2)$ and $(u^*_i)_{i=1,..,p}$
its dual basis. The element $\delta_U$ of $U o U^!$ verifies
${\delta}^{n+1}=0$.  Smirnov has defined a tensor product
$\otimes'$ in the category of operads such that there exists an
embedding $U o V\rightarrow U\otimes' V$ where $U$ and $V$ are
linear operads. We denote by $\delta'_U$ the image of $\delta$ by
the embedding $U o U^!\rightarrow U\otimes' U^!$, and by
$(U\otimes' U^!,\delta'_U)$ the resulting $n$-complex that we call
the $n$-Koszul complex. For $n=2$. This complex is the Koszul
complex defined by Ginzburg and Kapranov.

\bigskip

{\bf Representations of quadratic categories.}

\bigskip

{\bf Definition 1.}

Let $U=L[V]/C$ be a quadratic algebra. We denote by $Aut(U)$. The
automorphisms group of $U$. An element $h$ of $Aut(U)$ is defined
by an automorphism $h_1$ of $V$ such that $(h_1\otimes h_1)(C)=C$.
Consider an algebraic group $H$, a quadratic representation of
$H$, in $U$ is a morphism    $\rho:H\longrightarrow Aut(U)$.

\bigskip

{\bf Proposition 2.}

{\it Let $U$ be a quadratic algebra, and $H$ an algebraic group,
the set of quadratic representations $H\rightarrow Aut(U)$  is a
quadratic category.}

\bigskip

{\bf Proof.}

Let $\rho_1:H\longrightarrow Aut(U_1)$, and
$\rho_2:H\longrightarrow Aut(U_2)$ be two quadratic
representations, $\rho_1 o\rho_2$ is the representation defined
for each $h\in H$ by $(\rho_1 o \rho_2)(h):U_1 o
U_2\longrightarrow U_1 o U_2$ to be the morphism $\rho_1(h) o
\rho_2(h)$.

We define $(\rho_1\bullet\rho_2)(h)$ to be the automorphism of
$\rho_1(h)\bullet\rho_2(h)$ of $U_1\bullet U_2$.

We define by $(\rho_1^!)(h)$ is the automorphism $(\rho_1(h))^!$
of $U_1^!$.

The category $C_H$ of quadratic representations endowed with the
tensor products $o$ and $\bullet$ that we have just defined is a
quadratic category.

The respective neutral elements $I^H_o$ and $I^H_{\bullet}$ for
$o$ and $\bullet$ are the respective  trivial representations
$H\rightarrow Aut(I_o)$ and $H\rightarrow Aut(I_{\bullet})$.

\bigskip

Let $V$ be a vector space, $V$ is a quadratic algebra endowed with
the zero product. In this case $C=V^{\otimes^2}$. Thus the
category of representations of $H$ is embedded in the category of
quadratic representations of $H$.

\bigskip

{\bf Definition 3.}

A morphism of the category $C_H$ of quadratic representations is
defined by a family of automorphisms $e_U$ of $U$, where
$\rho:H\longrightarrow Aut(U)$ is an object of $C_H$, $U=T(U_1)/C$
such that $e_U$ commutes  with the action of $H$ on $U_1$. We
suppose that:

$e_{U\bullet V}=e_U\bullet e_V$

$e_{U o V}=e_U o e_V$

$e_{U^!}={e_U}^!$.

We denote $Aut(C_H)$ the group of automorphisms of $C_H$.

\bigskip

{\bf Proposition 4.}

{\it The natural embedding $H\rightarrow Aut(C_H)$ is an
isomorphism.}

\bigskip

{\bf Proof.}

The proof follows from the corresponding result for tensor
categories since the category of quadratic representations
contains the category of representations $\bullet$

\bigskip

{\bf Theorem 5.}

{\it Let $(C,o,\bullet,!)$ be a quadratic braided  abelian
category such that $End(I_{\bullet})$ is the ground field $L$, $P$
the category of quadratic algebras defined over the field $L$.
Suppose that there exists an exact faithful functor
$F\longrightarrow P$,
 then $(C,o,\bullet,!)$ is equivalent to the category of quadratic
 representations of an affine group scheme.}

\bigskip

{\bf Proof.}

Let $P'$ be the category whose objects are couples $(T,N)$ where
$T$ is a finite dimensional vector space, and $N$ a subspace of
$T\otimes T$. A morphism $h:(N,T)\longrightarrow (N',T')$ is a
linear map $h:N\longrightarrow N'$ such that $(h\otimes
h)(T)\subset T'$ The category $P'$ is endowed with the tensor
product defined by $(N,T)\otimes (N',T')=(N\otimes
N',t_{23}(T\otimes T'))$. This defines the morphism: $(h,h\otimes
h):(N,T)\longrightarrow (N',T')$. The functor $D:P\rightarrow P'$,
$U=T(U_1)/C\longrightarrow (U_1,C)$ is fully faithful, and it is a
tensor functor when $P$ is endowed with the tensor product
$\bullet_P$ defines by $T(U_1)/C\bullet_P T(U'_1)/C'=T(U_1\otimes
U'_1)t_{23}(C\otimes C')$.

Let $(C,o,\bullet,!)$ be a quadratic braided abelian category
endowed with an exact faithful quadratic functor $F:C\rightarrow
P$, then $D\circ F:C\rightarrow P'$ is exact and faithful. In
Deligne Milne, it is shown that $C$ is equivalent to the category
of comodules over a coalgebra ${\cal H}$. Since the object of $P'$
are finite dimensional spaces as in Deligne Milne, we can endow
${\cal H}$ with a structure of a bialgebra since every morphism
between two quadratic functors is an isomorphism, and deduce that
$(C,\bullet,o,!)$ is equivalent to the category or quadratic
representations of the affine group scheme $H$ defined by ${\cal
H}$ $\bullet$

\bigskip

 In the proof of theorem 5, we can define on $P'$ the following
 tensor structure: Let $(U,T)$, and $(U',T')$ be two elements of
 $P'$, $(U,T) o (U',T')=(U\otimes U',t_{23}(U\otimes T'+ T\otimes
 U'))$. This defines on the coalgebra ${\cal H}$ another algebra
 product.

 \bigskip

{\bf Definition 6.}

A neutral braided Tannakian quadratic category is a quadratic
category $(C,o,\bullet,!)$ such that there exists an exact
faithful $L$-linear functor $F:C\longrightarrow P$, where $P$ is
the category of quadratic algebras.   The functor $F$ is called a
fibre functor. If $L'$ is a $L$-algebra, a $L'$-fiber functor is a
fiber functor $F$ such that $D o F$ (see proof of theorem 5 for
the definition of D) takes its values in the category of
$L'$-modules.

\bigskip

{\bf Theorem 7.}

{\it Let $(C,o,\bullet,!)$ be a Tannakian quadratic category, $L$
a field and $Aff_L$ the category of affine  schemes defined over
$L$. We denote by $Fib(C)$ the category of fiber functors of $C$.
The functor $Fib(C)\rightarrow Aff_L$, $F\rightarrow H$ is a gerbe
over $Aff_L$ The fiber of $Spec(L')$ where $L'$ is a $L$-algebra
is the $L'$-valued fiber functor.}

\bigskip

{\bf Proof.}

Let $F':C(o,\bullet,!)\rightarrow P$ be another fiber
$L'$-functor, we have to show that $Hom(F',F')$ is representable
by a $H$-torsor over $Spec(L')$ where $H$ represents $F$. The
composition of a fiber functor and the functor $D$ defined in the
proof of theorem $5$ defines fiber a functor $C$ to the category
of $L'$-module. We   apply the corresponding result for Tanakian
category to $(F_0(C),\bullet)$.

\bigskip

Let $N$ be a projective scheme. we can realize $N$ as a projective
spectrum of a quadratic algebra $U_N$. We define $T(U_N)$ to be
the rigid quadratic generated by $T(U_N)$ and study the affine
group scheme which represents the automorphisms of $T(U_N)$.

\bigskip

{\bf Quadratic motives.}

\bigskip

 Let $L$ be a field, and $V_L$ the
category of projective schemes defined on $L$. A cohomological
functor $F:V\rightarrow (C,\otimes)$, where $(C,\otimes)$ is a
tensor abelian category,  which verifies standards properties that
verified cohomologies theories like Kunneth formula $F(U\times
U')=F(U)\otimes F(U')$. The category of motives
$F_0:V_L\rightarrow N_L$ is an initial objects in the category of
cohomological functors, that is for every cohomological functor
$F:V\rightarrow C$, there exists a functor of abelian categories
$F_C:N_L\rightarrow C$ such that $F=F_C\circ F_0$. The
construction of the category of motives is given in Saavedra, and
Deligne-Milne.

Let $U$ be a quadratic algebra. We denote by $Pro(U)$ the non
commutative projective scheme defined by $U$. This is the category
of $U$-graded modules up to the modules of finite length. The
category of coherent sheaves over every projective scheme is
equivalent to the non commutative projective scheme defined by a
quadratic algebra. This realization is not unique, that is the
category of coherent sheaves defined on a projective scheme can be
equivalent to the non commutative projective schemes defined by
two non isomorphic quadratic algebras.

We define the category of quadratic motives to be the category of
$L$-quadratic algebras this is an abelian category.
\bigskip

\bigskip

\centerline{\bf References.}

\bigskip

1. Beilinson, Ginzburg, Schechman. Koszul duality. Journal
Geometric Physics 1988.

2. Berger, $n$-Koszulity for non quadratic algebras. 239

3. Berger, Dubois-Violette, Wambst. Homogeneous algebras.
math.QA/020335

 4. Deligne, P. Milne, J., Ogus, Shih. Lectures notes 900. 1982

5. Ginzburg, V. Kapranov, M. Koszul duality for operads. Duke
Math. Journal 1994

6. Manin, Y. Quantum goups and non commutative geometry. Centre de
Recherches Mathematiques. Montreal.

7. Priddy, Koszul resolutions. Trans. American. Math. Society 152

\end{document}